\documentclass{article}
\usepackage{amsfonts}
\usepackage{amsmath}

\input{tcilatex}

\begin{document}

\title{Powers of doubly-affine integer square matrices with one non-zero
eigenvalue}
\author{Peter Loly, Ian Cameron and Adam Rogers \\
Department of Physics and Astronomy, University of Manitoba,\\
Winnipeg, Manitoba, Canada R3T 2N2}
\maketitle

\begin{abstract}
When doubly-affine matrices such as Latin and magic squares with a single
non-zero eigenvalue (1EV) are powered up they become constant matrices after
a few steps. The process of compounding squares of orders $m$ and $n$ can
then be used to generate an infinite series of 1EV squares of orders $mn$.

The Cayley-Hamilton theorem is used to understand this 1EV property, where
their characteristic polynomials: $x^{n}=\mathit{L}_{n}x^{n-1}$, have just
two terms, with $n$\ the order of the square and $\mathit{L}_{n}$\ the
row-column sum.
\end{abstract}

\section{Introduction}

We consider a diagonal Latin square [DLS] $\mathbf{A}=sud4a$\ from Cameron,
Rogers and Loly \cite{DMPS} [DMPS] of four symbols (here $1,2,3,4$ in
arithmetic progression) which has just one non-zero eigenvalue (1EV), and
then its matrix square $\mathbf{A}^{2}$ and cube $\mathbf{A}^{3}$:

\begin{equation}
\begin{bmatrix}
1 & 2 & 3 & 4 \\ 
3 & 4 & 1 & 2 \\ 
4 & 3 & 2 & 1 \\ 
2 & 1 & 4 & 3%
\end{bmatrix}%
,%
\begin{bmatrix}
27 & 23 & 27 & 23 \\ 
23 & 27 & 23 & 27 \\ 
23 & 27 & 23 & 27 \\ 
27 & 23 & 27 & 23%
\end{bmatrix}%
,%
\begin{bmatrix}
250 & 250 & 250 & 250 \\ 
250 & 250 & 250 & 250 \\ 
250 & 250 & 250 & 250 \\ 
250 & 250 & 250 & 250%
\end{bmatrix}%
\end{equation}

The intermediate square of two symbols still has rows ($R$), columns ($C$)
and diagonals (main $d1$, dexter $d2$) (RCD) with a\ common linesum ($100$),
but the final stage is a constant matrix with identical elements, with all
higher powers then constant matrices. We described $sud4a$ as a mini-Sudoku
in DMPS\ because of its design as a miniature version of order 2 Latin
squares as subsquares in a $4\times 4$ frame, contrasted with the well known
Sudoku solutions of a $9\times 9$ puzzle of a larger Latin square of 9
symbols ($1,2,...9$)\ with one of each in every 3-by-3 subsquare. This
result in ($1$) characterizes the present study and the sequence is quite
dramatic when each unique element is designated by a colour, declining from
four colours, through two, to a uniform single colour.

This study was prompted by a note to an online group of magic square
bloggers c. 27 September 2016 by Miguel Amela [AMA] , in \textit{Powers of
Associative Magic Squares }\cite{Miguel}, which\textit{\ }renewed interest
w.r.t \ "powering-up matrixwise" associative\textbf{\ }magic squares (also
called regular or symmetric in much of the earlier literature - see Loly,
Cameron, Trump and Schindel \cite{LAA2009}\ [LAA], which also includes
references to earlier literature). This is quite different from the issue of
multimagic squares which remain DDA when all their elements are squared,
cubed, etc. see Boyer \cite{Boyer} and Derksen, Eggermont and Van den Essen 
\cite{Derksen}.

In particular Miguel Amela was concerned to show counterexamples to
Proposition 6.1 by Cook, Bacon \ and Hillman \cite{CBH2010} [CBH]:

"Let $\mathbf{M}$\ be any associative $p\times p$\ magic square with magic
constant $m$. Then $\mathbf{M}^{2n+1}$\ is magic with magic constant $%
m^{2n+1}$, and $\mathbf{M}^{2n}$\ is semimagic."

We will show counterexamples, in part because we find examples that become
constant matrices at orders $n=4,5,8,16,...$.

CBH stated that the proof of their proposition " involves matrix algebra and
an analysis of the eigenvalues ...", but details were not given. An
alternation for odd versus even powers of associative magic squares was also
claimed. We will illuminate this claim with a variety of examples, while
extending beyond the specific constraints of associative magic squares, i.e.
those magic squares with antipodal pair sums of $1+n^{2}$ (or Latins with
pair sums $1+n$).

While CBH were concerned with the "magicness" of the powered-up squares, the
present work evolves from Amela's later observation [30 September 2016]
linking some counterexamples to certain order 4 associative magic squares
having the property of just one non-zero eigenvalue, "1EV" - again see our
LAA.

The corresponding characteristic polynomials are obtained by evaluating: $%
\det |\mathbf{A}-x\mathbf{I}_{n}|=0$, where $\mathbf{I}_{n}$ is the $n$-by-$%
n $ identity matrix. For the first three powers ($p=1,2,3$) of $sud4a$\ the
characteristic polynomials [CharPoly] are given in turn with their
corresponding eigenvalues ($\lambda _{i}$), singular values (SVs: $\sigma
_{i}$), \ and matrix rank $r$ (equal to the number of non-zero SVs), and sum
of the 4th powers of the second to last SVs, $R$.

\begin{center}
$%
\begin{tabular}{|l|l|l|l|l|l|}
\hline
$p$ & CharPoly & $\lambda _{i}$ & $\sigma _{i}$ & $r$ & $R$ \\ \hline
$1$ & $x^{4}=10x^{3}$ & $10,0,0,0$ & $10,4,2,0$ & $3$ & $272$ \\ \hline
$2$ & $x^{4}=100x^{3}$ & $100,0,0,0$ & $100,8,0,0$ & $2$ & $4096$ \\ \hline
$3$ & $x^{4}=1000x^{3}$ & $1000,0,0,0$ & $1000,0,0,0$ & $1$ & $0$ \\ \hline
\end{tabular}%
$

Table 1 - $sud4a$ - RCD linesum $(s_{4}=10)^{p}$. The sole eigenvalue and
leading singular value are equal to RCD (and the 1EV), $\lambda
_{1}=(10)^{p} $.

\textbf{Note the vanishing R-index which is the sum of the 4th powers of all
but the lead SV.}
\end{center}

In DMPS we introduced "singular value clans", which are characterised by a
particular set of singular values, as a powerful way to classify different
magic and Latin squares of the same order, and the $R$-index:

\begin{equation}
R=\dsum\limits_{i=2}^{n}\sigma _{i}^{4}\text{ }
\end{equation}%
which is an integer for all integer squares - see DMPS.

The rest of this paper continues with some preliminary matrix spectral
issues, before completing the order 4 Latin square issues, then moving on to
the parallel situation for magic squares, and after that the results of
compounding squares of orders $m$ and $n$, as an extension of Chan and Loly 
\cite{chan} in order to generate an infinite series of 1EV squares of orders 
$mn$, which the present authors\ know \cite{RCL} [RCL] will also have just
1EV.

The present\textbf{\ }approach builds on our previous 2004-7 studies of
magic square spectra in \cite{LAA2009} [LAA] which revealed several cases
with 1EV and included a discussion of non-diagonable matrices as well as the
insight gained by looking at Jordan Normal Forms, together with our more
recent 2013 extension which also included Latin squares in DMPS. For $sud4a$
these all have Mattingly's \cite{Mattingly}\ multiplicity ($\mu $)\ of zero
eigenvalues (here $\mu =3)$, originally for magic squares.

The number of SVs (and thus the matrix rank) decrease from $3$ through $2$\
to $1$ in this powering. We will later use the Cayley-Hamilton theorem \cite%
{HJ} which states that the characteristic polynomial\ is also satisfied by
the matrix itself, to understand this 1EV property. These 1EV squares
provide hitherto unexpected examples which become constant on matrix
powering.

Integer Latin and magic squares each have characteristic row and column sums
which are therefore doubly-affine [DA]. In addition magic squares have the
same sums for both major diagonals. All of our early DA squares can be
scaled up by constructing larger compound (or composite) DA squares of
multiplicative order following an ancient strategy extended by Chan and Loly
[CL]. Circa 2008 we realized that this preserved the 1EV property. Details
and more varieties (not needed here) are included in a companion paper [RCL].

\section{Terminology for magic and Latin squares}

Our interest in 1EV matrices began c. 2004 with magic squares as reported in
a conference in 2007 \cite{IWMS} [IWMS] and published in 2009 [LAA], but we
did not envisage this powering issue, even though we were aware of prior
interest in powering up magic squares as indicated in section 3.5 of that
work. However in September 2016 Miguel Amela \cite{Miguel} [AMA] rekindled
this by drawing our attention to squaring 8 associative (or symmetric or
regular) \cite{Andrews} order 4 magic squares which we had previously
identified as having 1EV \cite{IWMS},\cite{LAA2009}.

We are particularly concerned to be more precise with the use of both magic
and semimagic terms by grounding them in a long tradition that stretches
more than a century, e.g. Andrews\ classic work \cite{Andrews}, itself
largely a compilation of works by several authors who published in \textit{%
The Monist}.

We find it worth clarifying magic terminology since CBH appear to have
followed Stark \cite{Stark} who used "magic" when just the row and column
sums are equal, usually called semimagic, thus ignoring the two main
diagonals usually included in the definition of traditional magic squares.
Since having the same diagonal sums is a feature of the great majority of
studies on magic squares \cite{LEX}, we include the diagonals in our
definition, and thus use the familiar term of semi-magic squares for those
lacking one or both diagonals.

Here magic is reserved for "classic" magic squares of sequential integers
with the same RCD sums \cite{Andrews}, and semimagic whenever one or both
diagonals do not have the same linesum, so we leave "magic" and "semimagic"
for those which have non-sequential elements . We note that semimagic and
Latin squares are also doubly-affine [DA], while magic squares are diagonal
doubly-affine [DDA]. We use magic squares with elements from the consecutive
set $1,2,3,...n^{2}$, and our Latin squares have $n$ symbols in each row and
column whose elements run from $1,2,...n$, which are also the most common
historically. Then their linesums (a.k.a. magic constants, $\mathit{L}_{n}$,
in the abstract) for magic and Latin squares are respectively: 
\begin{equation}
S_{n}=\frac{n}{2}(1+n^{2}),\;s_{n}=\frac{n}{2}(1+n)
\end{equation}%
e.g. $S_{3}=15$, $S_{4}=34,...$ and $s_{2}=3,s_{3}=6,s_{4}=10,...$.,
replacing the generic $L_{n}$\ used in the Abstract. We note that the
elements used here are in arithmetic progression, but that there is some
interest in magic squares which are not, e.g. Fibonacci magic squares in CBH
which we examine near the end.

All magic squares have $8$ phases on rotation and reflection, and we have
previously noted in LAA some variations in spectral properties as a result
that turn out to be relevant here. Amongst the varieties of magic squares
several types have further properties due to additional constraints, e.g.
associated (or associative or regular) where all antipodal pairs have the
same sum, pandiagonal where all parallel broken diagonals also have the
magic linesum, ultramagic which are both associative and pandiagonal, etc.
see LAA. Further there are other squares of interest, particularly with Ben
Franklin's bent-diagonal squares which at least are semimagic due to their
half row and column sum property - see Schindel, Rempel and Loly \cite{PRSA}
[PRSA].

Our tests include orders $3,4,5,8$ and $16$ for magic squares,\textbf{\ }but
began with the Latin square $sud4a$\ which affords by far the simplest
example. Also we do not restrict ourselves to associative magic squares,
knowing already that other types in LAA exhibit 1EVs. Also our order four
1EV examples multiply up to order 16 by the basic compounding process
described in 2002 by Chan and Loly \cite{chan} [CL] which preserves the 1EV
property.

For the 880 order 4 magic squares we use the indexing common to the book by
Benson and Jacoby \cite{BJ}, as well as the web pages of Harvey Heinz \cite%
{LEX}. Note also that CBH use "regular" for sequential elements - we prefer
to reserve regular as an alternate to associative as per Andrews 1917, 2004 
\cite{Andrews}.

\section{Doubly-affine considerations}

It is helpful now to rephrase CBH's Proposition 1.1: "If $\mathbf{A}$\ and $%
\mathbf{B}$\ are $n\times n$\ semimagic matrices with magic constants $M$\
and $N$, respectively, then $\mathbf{AB}$\ and $\mathbf{BA}$\ are semimagic
with magic constant $MN$"\textbf{.}

A proof was given by CBH which remains valid if we edit it to apply to
doubly-affine square matrices as a new Proposition:

\textbf{PROPOSITION 1:} \textbf{If }$\mathbf{A}$\textbf{\ and }$\mathbf{B}$%
\textbf{\ are }$n\times n$\textbf{\ }DA\textbf{\ matrices with magic
constants }$M$\textbf{\ \ and }$N$\textbf{, respectively, then }$\mathbf{AB}$%
\textbf{\ and }$\mathbf{BA}$\textbf{\ are }DA\textbf{\ with magic constant }$%
MN$\textbf{."}

\subsection{Multiplying by constant square matrices}

CBH also noted that: "If $\mathbf{A}$ is magic and if some power of $\mathbf{%
A}$ is constant, then since a constant matrix times a magic square is
constant, it follows that all higher powers are magic."

Our Proposition 1 now takes care of multiplying by a constant square matrix.
Such a constant matrix can be represented by $c\mathbf{E}_{n}$\ where $c$\
is a multiplicative factor and $\mathbf{E}_{n}$\ is an $n\times n$\ unit
matrix consisting of all\ $1$'s.

In Tables below we will characterize the nature of the DA integer squares by
their "Type", \ i.e. magic or diagonal DA (DDA), and semimagic or simply DA.

Also the elements do not need to be integers for either the Latin or magic
squares, although that is usually the case. Both Latin squares and magic
squares are characterised by constant row and column linesums, and even
without the same linesum(s) for both diagonals, are both DA, hence its use
in the title. All magic squares are characterised by having the same linesum
on all rows, colums and both diagonals.

\section{Characteristic Polynomials and the Cayley-Hamilton Theorem}

In Table 1 we saw that $sud4a$ has rank $3$, with characteristic polynomial:

\begin{equation}
x^{4}-10x^{3}=0\text{, or }x^{4}=10x^{3}.
\end{equation}%
We will build more evidence later for the role of a single non-zero
eigenvalue (1EV) but for now turn to looking at the characteristic
polynomials for an understanding of the powering of 1EV magical squares to
constancy - see \cite{LAA2009}. For $sud4a$ in (1) has solutions $x=10$, the
linesum eigenvalue, and the triply degenerate $x=0$. So according to the
Cayley-Hamilton theorem \cite{HJ} one can substitute the matrix $\mathbf{A}$
for for $x$ in the characteristic equation (if $\mathbf{0}$ is the null
matrix and $\mathbf{I}$ is the identity matrix of the same\ order, here $n=4$%
, with $s_{1}=\lambda _{1}$):%
\begin{equation}
\mathbf{A}^{4}-10\mathbf{A}^{3}=0\text{, or }\mathbf{A}^{4}=10\mathbf{A}^{3}.
\end{equation}%
Which means that \textbf{both} $\mathbf{A}^{3}$\ AND $\mathbf{A}^{4}$\ are
constant matrices, and then multiplying both sides successively by $\mathbf{A%
}$, for $\mathbf{A}^{5}=10\mathbf{A}^{4}$, and so on for constancy at all
higher powers. We see that all powers higher than $p=2$ are \ constant
matrices with just one non-zero eigenvalue ($\mu =3$).

We conclude more generally that Cayley-Hamilton applied to 1EV
characteristic polynomials amounts to:

\textbf{PROPOSITION 2: Doubly-affine squares with 1EV power up }$\mathbf{A}$%
\textbf{\ to constant matrices.}

Which we will illustrate by examples where constancy occurs at either the
cube or quartic power.

\subsection{Singular Values (SVs) and Gramian matrix product}

The source of the matrix singular values (SVs), the\textbf{\ }Gramian matrix
product used to obtain the squares of the SVs from the product of a matrix
and its transpose, either\textbf{\ }$\mathbf{AA}^{T}$ or $\mathbf{A}^{T}%
\mathbf{A}$, is quite similar to the squared product of a matrix with
itself, with the same row and column sums and for $sud4a$ is:

\begin{equation}
\begin{bmatrix}
30 & 22 & 20 & 28 \\ 
22 & 30 & 28 & 20 \\ 
20 & 28 & 30 & 22 \\ 
28 & 20 & 22 & 30%
\end{bmatrix}%
\end{equation}%
while differing in the information that it yields. Those SVs lead us to the
entropic measures introduced by Newton and De Salvo \cite{NDS} [NDS] in 2010
for their Sudoku study, which we extended to a wider set of Latin squares as
well as magic squares in DMPS, and used here for comparing successive powers.

\section{Measures of the diversity of matrix elements}

Our results will be presented in terms of the "Compression" entropic measure
($C$) introduced for Sudoku solutions following NDS, which offers an insight
into the increase in uniformity upon powering up through squaring, cubing,
etc. through a measure based on the singular values of matrices.
Specifically one takes the Shannon entropy, $H$, which is a function of the
singular values normalized by their sum, $\widehat{\sigma }_{i}$, and then
used to obtain the percentage Compression $C$ in NDS:

\begin{equation}
H=-\Sigma _{i}\widehat{\sigma }_{i}\ln (\widehat{\sigma }_{i})\text{ , \ }%
C=(1-H/\ln (n))\ast 100\%
\end{equation}

An alternative "$Spread$" measure of the elements divides the difference
between max and min elements of matrices by their average value:%
\begin{equation}
Spread=n(Max[f_{i}]-Min[f_{i}])/L(n),
\end{equation}%
where $L(n)$ is the row-column linesum eigenvalue. For $sud4a$:

\begin{center}
$%
\begin{tabular}{|l|l|l|l|l|}
\hline
$p$ & $r$ & $C\%$ & $Spread$ & Type \\ \hline
$1$ & $3$ & $35.0603$ & $1.2$ & diagonal Latin (DLS) ($R=272$) \\ \hline
$2$ & $2$ & $80.9527$ & $0.16$ & DDA (\& associative) \\ \hline
$3$ & $1$ & $100$ & $0$ & constant; $250\mathbf{E}_{4}$ \\ \hline
\end{tabular}%
$

Table 2 - $sud4a$ - $p$ is the matrix power, $r$ the matrix rank, $C\%$ the\
Compression (here to 4 decimal places). \textbf{The spectra listed in Table
1 are not repeated here.}
\end{center}

\ In our tables $C\%=100$ means that every element of the cubed and higher
powers of $sud4a$ has all elements the same. Increasing powers have
decreasing $Spread$, although we prefer the increasing NDS Compression
measure.

\section{A note on Jordan Normal Form}

The 1EV squares that produce the constant matrix are non-diagonalizable [see
LAA], as revealed by the\textbf{\ }algebraic and geometric multiplicities of
the eigenvalues. In these situations a decomposition of the
non-diagonalizable matrix $\mathbf{A}$ is still possible via the Jordan
normal form, $\mathbf{A}=\mathbf{PJP}^{-1}$ where $\mathbf{P}$ is an
invertible matrix and $\mathbf{J}$ is the Jordan matrix \cite{HJ},\cite%
{Ortega}. We note that this form allows powers of the matrix to be computed
easily $\mathbf{A}^{n}=\mathbf{PJ}^{n}\mathbf{P}^{-1}$ , since we have $%
\mathbf{PP}^{-1}=\mathbf{I}$\textbf{\ }, the identity matrix.

The Jordan form of the square is useful since it is composed of $q$\ blocks
each of rank $r_{i},i=1..q$\textbf{,} whose sum corresponds to the rank of
the square and the maximum blocksize $k$ indicates the exponent required to
reduce our magic square to a constant matrix.

The Jordan normal form (often called the Jordan canonical form) is an upper
triangular matrix which has each non-zero off-diagonal entry equal to $1$,
immediately above the main diagonal (on the superdiagonal), and with
identical diagonal entries to the left and below them to form blocks, e.g.
the 3-by-3, 2-by-2 and 1-by-1 blocks:

\begin{center}
\begin{equation}
\begin{bmatrix}
\lambda _{1} & 1 & 0 &  &  &  &  &  \\ 
0 & \lambda _{1} & 1 &  &  &  &  &  \\ 
0 & 0 & \lambda _{1} &  &  &  &  &  \\ 
&  &  & \lambda _{2} & 1 &  &  &  \\ 
&  &  & 0 & \lambda _{2} &  &  &  \\ 
&  &  &  &  & 0 &  &  \\ 
&  &  &  &  &  & \lambda _{q-1} &  \\ 
&  &  &  &  &  &  & \lambda _{q}%
\end{bmatrix}%
.
\end{equation}
\end{center}

For a diagonalizable matrix the Jordan form simply has eigenvalues along the
main diagonal, with $1^{\prime }s$ above the diagonals in each block of
dimension greater than one and where all non-displayed entries are zero \cite%
{Olver}.

\section{Nilpotent Square Matrices with 1EV}

A nilpotent matrix is a square matrix $\mathbf{N}$ such that $\mathbf{N}%
^{q}=0$ for some positive integer $q$, where its smallest value is sometimes
called the degree or index of $\mathbf{N}$.

Let $\mathbf{Z}$ be a Latin or magic square of order $n$ with magic constant
linesum $l$ and $n-1$ zero eigenvalues, then we can write $\mathbf{Z}$\ in
terms of $\mathbf{N}$, a nilpotent matrix of order $n$ and index $k\leq n$
with zero line sum and the square matrix of order $n$ of all ones $\mathbf{E}%
_{n}$, in the following fashion:%
\begin{equation}
\mathbf{Z}=\mathbf{N}+(\frac{l}{n})\mathbf{E}_{n}\text{, or }\mathbf{N}=%
\mathbf{Z}-(\frac{l}{n})\mathbf{E}_{n}\text{.}
\end{equation}%
The square $\mathbf{N}$ will be nilpotent since it will have all zero
eigenvalues as required. If we raise $\mathbf{Z}$ to the $k$th power we will
end up with a polynomial in $\mathbf{N}$ and $\mathbf{E}_{n}$ composed of
terms $\mathbf{N}^{k}$, $(\frac{l}{n})^{k}\mathbf{E}_{n}^{k}$, and
cross-terms in $(\frac{l}{n})^{p}\mathbf{N}^{k-p}\mathbf{E}_{n}^{p}$, with
appropriate binomial coefficients. \ The first term $\mathbf{N}^{k}$
vanishes since we have index $k$, the crossterms will all vanish because
they involve multiplying $\mathbf{N}$ by $\mathbf{E}_{n}$ which simply gives
zero elements becauhe effect of $\mathbf{E}_{n}$ is to sum the rows/columns
of $\mathbf{N}$ that has zero linesum. Since $\mathbf{E}_{n}^{k}=n_{n}^{k-1}%
\mathbf{E}_{n}$, we end up for our magic square that $\mathbf{Z}^{k}=(\frac{%
l^{k}}{n})\mathbf{E}_{n}$, a constant matrix.

Also for a nilpotent matrix $\mathbf{N}$\ of index $k$ the Jordan form, $%
\mathbf{J}_{k}$, has just $1^{\prime }s$ above the main diagonal, from which
it follows that to reduce the Jordan form to a zero matrix we need to raise
it successively to the $k$th power while watching the one's migrate out of
the resulting products.

For the 1EV order 4 DA squares of interest here the Jordan form is:

\begin{equation}
\mathbf{J}_{4}=%
\begin{bmatrix}
0 & 1 & 0 &  \\ 
0 & 0 & 1 &  \\ 
0 & 0 & 0 &  \\ 
&  &  & \lambda%
\end{bmatrix}%
\text{, with }\lambda =l\text{, the linesum eigenvalue}.
\end{equation}

The largest blocksize is 3 which has rank 2, so the rank of this "4-square"
is $r_{4}=2+1=3$ . If $L_{4}$ is the DA linesum the diagonal elements for $%
\mathbf{J}_{4}^{3}$ are $\{0,0,0,L_{4}^{3}\}$, and all superdiagonal
elements are zero. This form is characteristic of all $n=4$ 1EV magic
squares, and a $\mathbf{J}_{5}$ example for $n=5$ will be shown later.

\section{Jordan forms for powers of $sud4a$ ($\protect\mu =3$)}

For $sud4a$ and its powers in the Introduction the second power has all
element pairs symmetric about the centre with the same sum, there $50$, and
can be classified as associative (or regular) using terminology from\ the
literature of magic squares \cite{LAA2009}, with higher powers trivially so.
However the initial $sud4a$ is clearly not associative.

We now show their Jordan Forms \cite{Ortega},\cite{HJ} in terms of blocks
along the main diagonal:

\begin{equation}
\begin{bmatrix}
0 & 1 & 0 &  \\ 
0 & 0 & 1 &  \\ 
0 & 0 & 0 &  \\ 
&  &  & 10%
\end{bmatrix}%
,%
\begin{bmatrix}
0 & 0 & 0 &  \\ 
0 & 0 & 1 &  \\ 
0 & 0 & 0 &  \\ 
&  &  & 100%
\end{bmatrix}%
,%
\begin{bmatrix}
0 & 0 & 0 & 0 \\ 
0 & 0 & 0 &  \\ 
0 & 0 & 0 &  \\ 
&  &  & 1000%
\end{bmatrix}%
\end{equation}%
Here the off-diagonal "$1$'s" move out for each rise in power, with their
matrix ranks first from an order $3$ block of rank $2$ and the diagonal "$10$%
" for rank $r=2+1=3$, next an order $2$ block and the diagonal "$100$" for
rank $r=1+1=2$, and then final rank $r=0+1=1$, after which the higher powers
are all constant matrices.

\section{lat4a - 3 non-zero eigenvalues [3EV] ($\protect\mu =1$)}

There are many singular squares of order 4 and higher with three non-zero
eigenvalues, the abbreviation 3EV will be useful. Here we compare $sud4a\ $%
with $lat4a$, from which we obtained the former by moving the second row to
the last. Now the first three powers of $lat4a$\ are:

\begin{center}
\begin{equation}
lat4a=%
\begin{bmatrix}
1 & 2 & 3 & 4 \\ 
2 & 1 & 4 & 3 \\ 
3 & 4 & 1 & 2 \\ 
4 & 3 & 2 & 1%
\end{bmatrix}%
,%
\begin{bmatrix}
30 & 28 & 22 & 20 \\ 
28 & 30 & 20 & 22 \\ 
22 & 20 & 30 & 28 \\ 
20 & 22 & 28 & 30%
\end{bmatrix}%
,%
\begin{bmatrix}
232 & 236 & 264 & 268 \\ 
236 & 232 & 268 & 264 \\ 
264 & 268 & 232 & 236 \\ 
268 & 264 & 236 & 232%
\end{bmatrix}%
,...
\end{equation}
\end{center}

None of these are associative, but all are antipodal about the center [or
symmetric about both main diagonals], and all still four symbol Latin. The
spectra of $lat4a$ and its first four powers are:

$p=1$:\textbf{\ }$\lambda _{i}=10,-4,-2,0$; $\sigma _{i}=10,4,2,0$;

$p=2$: $\lambda _{i}=$ $\sigma _{i}=100,16,4,0$;

$p=3$: $\lambda _{i}=10^{3},-64,-8,0$; $\sigma _{i}=$ $10^{3},64,8,0$;

$p=4$: $\lambda _{i}=\sigma _{i}=$ $10^{4},256,16,0$; ...

all with multiplicity $\mu =1$ of zero eigenvalues. Then other measures are:

\begin{center}
$%
\begin{tabular}{|l|l|l|l|l|l|}
\hline
$p$ & $r$ & $C\%$ & $Spread$ & Type & $R$ \\ \hline
$1$ & $3$ & $35.0603$ & $1.2$ & Latin $d1=4,d2=16$ & $272$ \\ \hline
$2$ & $3$ & $61.4828$ & $0.4$ & DA $d1=120,d2=80$ & $65,792$ \\ \hline
$3$ & $3$ & $80.5474$ & $0.144$ & DA $d1=928,d2=1072$ & $16,781,312$ \\ 
\hline
$4$ & $3$ & $90.7518$ & $0.0544$ & DA $d1=10272,d2=9728$ & $4,295,032,832$
\\ \hline
\end{tabular}%
\ \ \ \ \ \ \ \ \ \ \ \ $

Table 3 - $lat4a$
\end{center}

\bigskip In Table 3 note how $Compression$ gradually increases towards the
100\% of uniformity, while the $Spread$ decreases gradually towards
uniformity, but neither is ever achieved. Here the average of the diagonal
sums is the common RC linesum in Table 1\textbf{. }

The Jordan Normal form was already useful to us for non-diagonalizable $%
sud4a $, as shown in \cite{LAA2009} for magic squares, so the diagonal form
for $lat4a$ which is diagonalizable has eigenvalues $-4,-2,0,10$\ down the
main diagonal, $r=3$.

A discussion of eigenvalue and singular value spectra for low order Latin
and magic squares that is useful in the present context is found in \cite%
{DMPS}.

\section{The $n=3$\ associative magic square ($\protect\mu =0$)}

Aside from rotations and reflections there is just one third order magic
square, the ancient Loshu. Here $\mathbf{A}=loshu=\{\{4,9,2\},\{3,5,7\},%
\{8,1,6\}\}$, which is the top-bottom reflection of $\mathbf{A}_{1}$ of CBH.
We show it followed by its square, cube and quartic matrix powers ($p=1,2,3$%
):

\begin{equation}
\begin{bmatrix}
4 & 9 & 2 \\ 
3 & 5 & 7 \\ 
8 & 1 & 6%
\end{bmatrix}%
,%
\begin{bmatrix}
59 & 83 & 83 \\ 
83 & 59 & 83 \\ 
83 & 83 & 59%
\end{bmatrix}%
,%
\begin{bmatrix}
1149 & 1029 & 1197 \\ 
1173 & 1125 & 1027 \\ 
1053 & 1221 & 1101%
\end{bmatrix}%
,...
\end{equation}

Here the non-singular characteristic polynomials are respectively:

$p=1$:\textbf{\ }$\lambda _{i}=$ $15,2i\sqrt{6},-2i\sqrt{6}$; $\sigma
_{i}^{2}=225,48,12;R=2448$;

$p=2$: $\lambda _{i}=$ $225,-24,-24$; $\sigma _{i}=225,24,24$;\textbf{\ note
the degeneracy in both spectra}.\textbf{\ }Here\textbf{\ }$R=663,552$;

$p=3$: $\lambda _{i}=$ $3375,\sqrt{13824}I,-\sqrt{13824}I$; $\sigma
_{i}=3375,\sqrt{27648},\sqrt{6912}$;\textbf{\ }$R=812,187,648$.

While the matrix elements increase rapidly, notice how the relative spread
of the elements progressively decreases on powering up. These are tabled up
to sixth matrix power in the following table$\bigskip $

\begin{center}
$%
\begin{tabular}{|l|l|l|l|l|}
\hline
$p$ & $r$ & $C\%$ & $Spread$ & Type \\ \hline
$1$ & $3$ & $14.7017$ & $1.4$ & magic ($R=2448$) \\ \hline
$2$ & $3$ & $46.5804$ & $0.32$ & DA, $d1=177$ \\ \hline
$3$ & $3$ & $73.2057$ & $0.1707$ & DDA \\ \hline
$4$ & $3$ & $88.8869$ & $0.0341$ & DA, $d1=51777$ \\ \hline
$5$ & $3$ & $95.3843$ & $0.0182$ & DDA \\ \hline
$6$ & $3$ & $98.2995$ & $0.0036$ & DA, $d1=11362977$ \\ \hline
\end{tabular}%
$

Table 4 - $n=3$ Loshu magic square - the linesum at each power is\textbf{\ }$%
S_{3}=15^{p}$, i.e. $15,225,3375,50625,..$\textbf{\ }
\end{center}

In Table 4 the odd and even powers are in agreement with the alternation in
CBH's Propositions 2.1 for 3-by-3's, and their general Proposition 6.1, as
their "magic" alternates with "semimagic" as powers increase and the
Compression monotonically tends to 100\%, while the $Spread$ tends
progressively to zero.

In LAA we showed the Jordan Form of the $loshu$, which has just the
eigenvalues on the main diagonal, and also discussed the 8 phases of magic
squares using the Loshu as an example of changing eigenvalues from $15,2i%
\sqrt{6},-2i\sqrt{6}$ to $15,2\sqrt{6},-2\sqrt{6}$\ on rotation by 90
degrees\ in\ the first table on page 2667 of LAA.

\section{Prior applications and studies of matrix powering}

The previous 3EV cases of $lat4a$ and $loshu$ both illustrate the gradual
increase in Compression and reduction in Spread which are more
representative of larger general (not doubly-affine) matrices where the
iterative powering up of matrices towards constancy has proved useful, e.g.
for the largest eigenvalue of large matrices - Pitman \cite{Pitman}
reference to Von Mises, etc. with addition of the corresponding eigenvector
evolutions.

The Compression and Spread measures in Tables 2 and 3 are an indicator of
what happens in the cases of larger matrices of full rank, but singular
magic and Latin squares may coverge much more rapidly as shown by some
examples shown below. First we use our earlier results for the spectral
properties of the well known complete set of 880 fourth order magic squares
where we elucidated the 1EV phenomenon [LAA], and where we noted earlier
powering literature, before a brief discussion of 1EV fifth order magic
squares from the complete set of some 275 million. Higher orders are
impossible to survey in their entirety, but prompted by CBH we take a look
at order 8 in the context of Franklin squares [PRSA].

\section{n=4 magic squares - Dudeney Groups and SV clans}

From Dudeney 1917 \cite{Dudeney} we have an important classification of
Frenicle's 1699 count into the Dudeney Group classification listing of the
880 distinct order 4 magic squares. We refer to the Appendix of Benson and
Jacoby \cite{BJ} for a list of the 880 with their Dudeney Group designation.
In 2007-9 LAA noted [\^{\i}n that reference's table 3] that the 7 Dudeney
groups [I,II,III,IV,V,VI-P,VI-S] were singular with at least one zero
eigenvalue. Then in DMPS we gave a summary of order $n=4$ magic squares
which have 880 distinct members in 13 Dudeney groups and the then new 63 SV
clans which we introduced.

Further Dudeney Groups I (pandiagonal), II (semi-pandiagonal, semi-bent) and
III (associative/regular, semi-pandiagonal), each have 48 distinct members,
each with the same three singular value clans, $\alpha ,\beta ,\gamma $,
characterised by their singular values, as shown in DMPS in that reference's
Figure 2.

In the associatives (Group III) there are eight 1EV squares:

\begin{center}
$%
\begin{tabular}{|l|l|l|l|}
\hline
clan & $\sigma _{i}$ & Frenicle indices & $R$ \\ \hline
alpha & $34,8\sqrt{5},2\sqrt{5}$ & $290,360,790,803$ & $102,800$ \\ \hline
beta & $34,4\sqrt{17},2\sqrt{17}$ & $99,377,489,535$ & $78,608$ \\ \hline
gamma & $34,2\sqrt{5}\sqrt{13},4\sqrt{5}$ & e.g. $126$ (of forty total) & $%
74,000$ \\ \hline
\end{tabular}%
\ \ \ \ \ \ $

Table 5 - summary of Group III SV clans.
\end{center}

For comparison with CBH we begin with one which agrees with their
Proposition 6.1, choosing the 1EV $f360$ which is the transpose of AMA's
fourth order associative example. AMA noted that $8$\ Group III Associative
MSs with multiplicity $m=3$\ of zero eigenvalues [now identified as those
with the single linesum eigenvalue, 1EV] as described in LAA. AMA's\textbf{\ 
}counterexample to CBH's Prop. 6.1, although AMA only took this to the
second power and did not note the constancy which begins at power\ $p=3$. $%
f360$\textbf{\ }is also a column permutation of MATLAB's magic(4) which is
(29) in LAA, p.2668: col.1-\TEXTsymbol{>}col.3; col. 2-\TEXTsymbol{>}col.1;
col.3-\TEXTsymbol{>}col.4; col.4-\TEXTsymbol{>}col.2.]

\subsection{$f360$ clan alpha associative ($\protect\mu =3$) 1EV}

Now:

\begin{center}
\begin{equation}
\begin{bmatrix}
2 & 11 & 7 & 14 \\ 
13 & 8 & 12 & 1 \\ 
16 & 5 & 9 & 4 \\ 
3 & 10 & 6 & 15%
\end{bmatrix}%
,%
\begin{bmatrix}
301 & 285 & 293 & 277 \\ 
325 & 277 & 301 & 253 \\ 
253 & 301 & 277 & 325 \\ 
277 & 293 & 285 & 301%
\end{bmatrix}%
,9826\mathbf{E}_{4}
\end{equation}
\end{center}

The eigenvalues of these matrices are respectively: $34,0,0,0$; $1156,0,0,0$%
; and $39304,0,0,0$, note how this becomes a matrix of identical elements
after power $2$, as did $sud4a$, and SVs - clan alpha\textbf{\ }in Table 5.
The linesums are\textbf{\ }$S_{4}^{p}=34^{p}$, i.e. $34,1156,39304$. Their
characteristic polynomials have the form: $-S_{4}x^{3}+x^{4}=0$, or $%
x^{4}=S_{4}x^{3}$, where $S_{4}=34$, with eigenvalues $\lambda _{i}=34,0,0,0$%
, and similarly for the square and cube with 34 replaced by 1156 and 39304
respectively. The next table summarizes the results up to matrix power $p=3$:

\begin{center}
$%
\begin{tabular}{|l|l|l|l|l|l|}
\hline
$p$ & $r$ & $C\%$ & $Spread$ & Type & $R$ \\ \hline
$1$ & $3$ & $37.2284$ & $1.7647$ & magic & $102,800$ \\ \hline
$2$ & $2$ & $82.7039$ & $0.2491$ & DDA & $40,960,000$ \\ \hline
$3$ & $1$ & $100$ & $0$ & constant: \ $9826\mathbf{E}_{4}$ & $0$ \\ \hline
\end{tabular}%
\ \ \ \ \ \ $

Table 6 - $f360$ Group III, clan alpha - Note the vanishing R-index.
\end{center}

This does not agree with CBH's Proposition 6.1 in their alternation \
between "magic" and "semimagic", in fact the behaviour is totally different
because of the constancy and since for $p=2$ we find our DDA instead of
their "semimagic". Note how the rank descends from 3 through 2 to rank 1,
i.e. a matrix of equal elements, and remains there for higher powers.

N.B. CBH late in their paper in Proof of their Proposition 5.1 for an
unusual order 8 magic square also indicate this point, but did not link it
to the 1EV property. More on this later.

\subsection{$f299$ clan beta associative ($\protect\mu =3$) 1EV}

Consider also a second 1EV from Group III, but now in SV clan\ beta:

$f299=\{\{2,7,13,12\},\{16,9,3,6\},\{11,14,8,1\},\{5,4,10,15\}\}$.

Now the only differences from Table 6 are a slightly different elements on
squaring, and compressons: $C\%=31.5627$\ and $75.7272$, as well as a $%
Spread $ of $0.3460$\ for $p=2$.

\subsection{Associative ($\protect\mu =1$) 3EV magic square $f175$}

For comparison we now look at one of the forty 3EV Associative (Regular,
semi-pandiagonal) Group=III in LAA's ($31)$ which affords an example on a
par with CBH.\ Table 5 in LAA shows that it is a member of SV clan alpha in
Table 5. $\mathbf{A}=f175$ is shown next, followed by its matrix square and
cube:

\begin{center}
\begin{equation}
\begin{bmatrix}
{\small 1} & {\small 12} & {\small 8} & {\small 13} \\ 
{\small 14} & {\small 7} & {\small 11} & {\small 2} \\ 
{\small 15} & {\small 6} & {\small 10} & {\small 3} \\ 
{\small 4} & {\small 9} & {\small 5} & {\small 16}%
\end{bmatrix}%
{\small ,}%
\begin{bmatrix}
{\small 341} & {\small 261} & {\small 285} & {\small 269} \\ 
{\small 285} & {\small 301} & {\small 309} & {\small 261} \\ 
{\small 261} & {\small 309} & {\small 301} & {\small 285} \\ 
{\small 269} & {\small 285} & {\small 261} & {\small 341}%
\end{bmatrix}%
{\small ,}...
\end{equation}
\end{center}

For $f175$: $\lambda _{i}=$ $34,8,-8,0$; $\sigma _{i}^{2}=\{1156,320,20,0\}$
showing the singularity, while the results are otherwise qualitatively
similar to those in Table 4 for n=3 with $C\%$ rising from $37.2284$ at $p=1$%
\ to $99.7369$\ at $p=6$\ and $Spread$\ declining from $1.7647$\ to $0.0009$%
, with the Jordan Form for $f175$ just has the eigenvalues on the diagonal.

While it is not our goal to give a complete account of all $n=4$ magic
squares, it does seem worthwhile to note other 1EV cases which are not
associative magic squares and thus become constant under powering. These are
found in Dudeney's \cite{Dudeney}\ Group I Pandiagonals and his Group II
Semi-pandiagonal and Semi-bent [see Table 3 in LAA], of which we choose one
of the latter next.

\subsection{$f27$ Group II semi-pandiagonal, semi-bent rank oscillator ($%
1\leftarrow \protect\mu \rightarrow 3$) - not associative}

A question that must be addressed is whether associativity in magic squares
such as $f360$\ is necessary for the powering pattern given by CBH of
alternation of "magic" and "semimagic". In Dudeney's Group II
[semi-pandiagonal, semi-bent magic squares, LAA] there are 32 magic squares
with multiplicity of zero eigenvalues $\mu =1$, and 16 oscillating on
reflection or rotation [LAA\textbf{\ }Table 3] between $\mu =1$\ and $\mu =3$%
, e.g. from the list in Benson and Jacoby \cite{BJ} we take first $f27$ with 
$\mu =1$ with eigenvalues $\lambda _{i}=\{34,-8,8,0\}$ which is similar to $%
f175$ in behaviour under powering, with $f27$ having a diagonal Jordan Form
the same as for $f175$ above so not shown explicity. Here the powering
details (omitted) follow CBH's Proposition 6.1, thus extending its scope
beyond just the associative magic squares in Group II.

\subsubsection{$f27$\textbf{\ rotated pi/2}}

Now with a different result compared to $f27$ but with the pattern of
constancy by power $p=3$ as in Table 6 for $f360$, and with the exception
only of a slightly different Type for $p=2$ of DA with $d2=1092$ instead of
DDA. So Dudeney \cite{Dudeney}\ Group II "semi-bent, semi-pandiagonals" have
the same behaviour for the rotated phase as did the Group III associatives
for all their phases. Thus 1EV magic squares do not need to be associative!\ 

We note also that Group I Pandiagonals have another 16 phase oscillators\ -
see Table 3 in LAA for sixteen more 1EV magic squares.

\subsection{ Non-singular Group IX $f181$ ($\protect\mu =0$) - not
associative}

For the sake of comparison with previous singular case we choose $\mathbf{A}%
=f181=\{\{1,12,13,8\},\{16,9,4,5\},\{2,7,14,11\},\{15,6,3,10\}\}$, taken
from equation (39) in LAA. The characteristic polynomial is $%
(x-34)(x+8)(x^{2}-8x+24)=0$, and now a diagonal Jordan form holds the
eigenvalues. Now we tabulate the properties through to the fifth power:

\begin{center}
$%
\begin{tabular}{|l|l|l|l|l|}
\hline
$p$ & $r$ & $C\%$ & spread & Type \\ \hline
$1$ & $4$ & $26.7676$ & $1.76471$ & magic ($R=92,288$) \\ \hline
$2$ & $4$ & $68.2067$ & $0.33218$ & DA $d1=1236,d2=1124$ \\ \hline
$3$ & $4$ & $90.4131$ & $0.08060$ & DA $d1=38728,d2=39496$ \\ \hline
$4$ & $4$ & $97.1769$ & $0.02011$ & DA $d1=1339536,d2=1335056$ \\ \hline
$5$ & $4$ & $99.3030$ & $0.00441$ & DA $d1=45397024,d2=45449248$ \\ \hline
\end{tabular}%
\ \ \ \ \ \ \ $

Table 7 - $f181$ non-singular and not associative.
\end{center}

Because $f181$\ is not an associative magic square this was not expected to
agree with CBH's Proposition 6.1 in their alternation \ between "magic" and
"semimagic".

\section{$n=5,$ $\protect\mu =4$ 1EV associative magic square}

Although there are no singular order 5 Latin squares [DMPS], our LAA paper
with Walter Trump already provides an order 5 magic square with just 1EV
(one just four cases),

\begin{center}
$\mathbf{A}=laa44=%
\begin{bmatrix}
2 & 11 & 21 & 23 & 8 \\ 
16 & 14 & 7 & 6 & 22 \\ 
25 & 17 & 13 & 9 & 1 \\ 
4 & 20 & 19 & 12 & 10 \\ 
18 & 3 & 5 & 15 & 24%
\end{bmatrix}%
,$
\end{center}

for which the characteristic polynomial for $p=1$ gives by Cayley-Hamilton $%
\mathbf{0=A}^{5}-65\mathbf{A}^{4}$ or $\mathbf{A}^{5}=65\mathbf{A}^{4}$, so
that by contrast with the earlier $n=4$ 1EV cases of $sud4a$ and $f27$ the
constancy occurs first at $p=4$ instead of $p=3$ (as found previously):

\begin{center}
$%
\begin{tabular}{|l|l|l|l|l|l|}
\hline
$p$ & $r$ & $C\%$ & spread & Type & $R$ \\ \hline
$1$ & $4$ & $25.3638$ & $1.84615$ & magic & $706,000$ \\ \hline
$2$ & $3$ & $64.2240$ & $0.41894$ & DDA & $82,414,854,400$ \\ \hline
$3$ & $2$ & $92.6908$ & $0.10487$ & DDA & (very large) \\ \hline
$4$ & $1$ & $100$ & $0$ & $3570125\mathbf{E}_{5}$ & $0$ \\ \hline
\end{tabular}%
\ \ \ \ $

Table 8 - A 1EV $n=5$ magic square ($\mu =4$) with R,C linesums of $%
(S_{5}=65)^{p}$.

Note again the vanishing R-index for a 1EV, but now at $p=4$.
\end{center}

This shows its rank 4 dropping by one in each step until it is a constant
matrix at $p=4$, instead of constancy at $p=3$\ as for $f360$. This
associative magic square differs from CBH's Proposition 6.1. Now the powers
of the Jordan Forms see the off-diagonal $1$'s moving out until constancy is
reached in a similar fashion to (12), i.e. a single 4-by-4 block of rank 3
for the zero eigenvalues, whose rank reduces by 1 for each increase in
power, with powers of $65$ left on the lower right corner. This shows why
constancy is attained at $p=4$\ rather than the previous $p=3$\ examples.

\subsection{$n=5,$ ($\protect\mu =0$) Ultramagic non-singular magic square}

In LAA\ we examined $laa45$, an ultramagic\ square (associative and
pandiagonal). Here the results are otherwise qualitatively similar to those
in Table 4 for $n=3$ with $C\%$ rising $p$\ as $23.4758$ to $99.4352$\ at $%
p=5$\ and $Spread$\ declining from $1.84615$\ to $0.003$. This agrees with
CBH's 2010 Proposition $6.1$\ which predicts alternation between their
"magic", our DDA, and their "semimagic", our DA.

We also note that there are no singular $n=5$ Latin squares DMPS, and so no
1EVs [DMPS].

\section{Are there 1EV squares in order 6 and higher?}

From DMPS there are singular rank 4 order 6 Latin squares but we are not
aware of any 1EV examples. Moreover there are no $n=6$ associative magic
squares - see Trump \cite{Trump}.

For $n=6$ and higher there are too many magic squares to count \cite{Trump},
but some of these are discussed in our LAA.\ So special cases must be
considered, one of which we address next.

\subsection{n=8 Franklin bent diagonal squares - rank 3 ($\protect\mu =5$)}

Here we consider questions concerning the complete set of $n=8$ Franklin
squares which were first counted by Schindel, Rempel and Loly [SRL] in \cite%
{PRSA} [PRSA]. These are characterised by common half row and column sums
(and thus are at least semimagic) as well as a complete set of $4n$ bent
diagonals [PRSA]. Franklin's handful of squares also had another property
that all 2-by-2 quartets have the half row-column sum, and which is
considered a further constraint. SRL were also the first to find magic
Franklin squares since exactly 1/3 of their definitive count were also DA in
their rows and columns. Moreover in DMPS\textbf{\ }we stated on p. 2675 that
all 368,640 eighth order natural pandiagonal Franklin squares have three
non-zero eigenvalues and that they are singular with rank 3.

Our landmark count of the order 8 Franklin bent diagonal semimagic squares
did not consider their spectral properties [PRSA], but in DMPS\ we reported
that they are rank 3 and thus highly singular.

CBH drew attention to a particular magic square, $BF$, which they derived
from a bona fide Franklin square of Neal Abrahams on Suzanne Alejandre's web
site \cite{abrahams}\ by RC permutations of rows 1,2 and then rows 3,4. For
Abrahams square (not shown here): $\lambda _{i}=$ $%
\{260,-43.7128,11.7128,0,0,0,0,0\},$ $\sigma _{i}^{2}=\{67600,4[2(1365\pm 
\sqrt{1,755,705})$ for $21520.2$ and $319.758\},$ the latter for integer $%
R=463,223,040$ and $r=3$,, which is consistent with Table 3 in DMPS, and
which is neither magic nor 1EV: $d1=228,d2=292$.

\subsubsection{CBH's "Franklin" 8th order 1EV magic squares ($\protect\mu =7$%
)}

CBH claimed to study "BEN FRANKLIN'S 8X8 MAGIC SQUARE MATRIX, $BF$. Added to
the puzzle in retrospect is Proposition 5.1 of CBH: \textbf{"The Alejandre
Form of the Ben Franklin }$8$\textbf{-by-}$8$\textbf{\ Magic Square and all
of its }$k$\textbf{th integral powers are magic, with magic constant }$%
(260)^{k}$\textbf{"}.

Their proof noted that for\textbf{\ }$k>2$ the powered up matrix was a
constant matrix (all elements equal), but this was understated and they did
not appear to extend this to other magic squares, and then we found that
their's was not a bona fide Franklin bent-diagonal square! \ Unfortunately
this is not a true Franklin square [see PRSA] because while magic (in RCD),
it is not so in all the bent diagonals. CBH\textbf{\ }use a row permutation
of Abrahams' bona fide Franklin square, which is not magic in the diagonals%
\textbf{\ }\cite{abrahams}, involving swapping row 1 and 2 and then row 3
and 4. While a bona fide magic square [RCD], in their proof of their
Proposition 5.1 CBH state that powers of this magic square for the cube and
higher are constant (magic) matrices with magic constant $(260)^{k+1}$, $k>2$%
. After correcting the last cell in CBH in the first row from a duplicate 10
to 19:

\begin{center}
\begin{equation}
BF=%
\begin{bmatrix}
14 & 3 & 62 & 51 & 46 & 35 & 30 & 19 \\ 
52 & 61 & 4 & 13 & 20 & 29 & 36 & 45 \\ 
11 & 6 & 59 & 54 & 43 & 38 & 27 & 22 \\ 
53 & 60 & 5 & 12 & 21 & 28 & 37 & 44 \\ 
55 & 58 & 7 & 10 & 23 & 26 & 39 & 42 \\ 
9 & 8 & 57 & 56 & 41 & 40 & 25 & 24 \\ 
50 & 63 & 2 & 15 & 18 & 31 & 34 & 47 \\ 
16 & 1 & 64 & 49 & 48 & 33 & 32 & 17%
\end{bmatrix}%
\text{.}
\end{equation}
\end{center}

We note that the modified square has just a single non-zero eigenvalue, in
contrast with three in the rank 3 Franklin n=8 set \cite{DMPS},\cite{PRSA}.
It has characteristic polynomial and spectra:

$-260x^{7}+x^{8},\lambda _{i}=260,0,0,0,0,0,0,0$; with the same SVs as above
for Abrahams' square, where we note as an aside that $\sigma _{i}^{2}$\ are
not always integer, and again with $R=463,223,040$, which is identical with
Abrahams square above, i.e. both are in the same SV clan [DMPS].

This all suggests to us that many more 1EV squares may be expected by
applying similar row and/or column permutations to other bona fide Franklin
squares. The Jordan form of this remarkable magic square is simply an
enlargement to $n=8$ of the $n=4$ Jordan forms in (11) for $sud4a$, $\mathbf{%
J}_{4}$, \ etc. by the addition of 4 more rows and columns of zeroes, with
the magic eigenvalue as appropriate, now $260$, and so need not take up
display space! Moreover the powering is similar to Table 6 with $%
C\%=61.4849,97.872,0;Spread=1.9385,0.0152,0.0$, with linesums\textbf{\ }$%
(S_{8})^{p}=260^{p}$, i.e. $260,67600,..$.

\begin{center}
$%
\begin{tabular}{|l|l|l|l|l|}
\hline
$p$ & $r$ & $C\mathbf{\%}$ & spread & Type \\ \hline
$1$ & $3$ & 61.4849 & 1.93846 & magic ($R=463,223,040$) \\ \hline
$2$ & $2$ & 97.872 & 0.0151479 & DDA \\ \hline
$3$ & $1$ & 100 & 0 & constant $\mathbf{E}_{8}$ \\ \hline
\end{tabular}%
\ \ \ \ $

Table 9 - 1EV $n=8$\ magic square,\textbf{\ }$BF$\textbf{.}
\end{center}

N.B. For "the proof" of Proposition 6.1 CBH\ state that it involves matrix
algebra and analysis of the eigenvalues, so perhaps they found the single
eigenvalue property?

Nordgren \cite{Nordgren} recently showed that odd matrix powers of true
pandiagonal Franklin squares [PRSA] are also pandiagonal Franklin squares,
and confirmed the rank 3 property.

\section{Compounding sud4a to order 16 1EV DLS}

Compounding is an ancient idea (more than a thousand years) for generating
larger magic squares, initially used for multiplying up the non-singular $%
n=3 $ $loshu$ to order $9$. A paper with Wayne Chan \cite{chan} [CL] on
rediscovering the basic compound idea used here, included iterating a
pandiagonal and most-perfect $n=4$ with Euler's 1779 pandiagonal $n=7$ to
the highly compounded $n=12,544=4^{4}7^{2}$, also pandiagonal.

Now the 1EV property extends to compounded 1EV Latin and magic squares. Next
we apply compounding to our 1EV Latin square, $sud4a$, before considering
the magic square parallels \cite{RCL} [RCL]. The order of the squares
increases to multiples of $n=4$, i.e. $n=16,64,256,...$ of which the first
will be considered in detail here.

In 2004 we began work on an expansion of the compounding idea to include
Latin squares, and general magic squares, which lead to a formula for the
eigenvalues and singular values of compounded magic (and Latin) squares
which was included in our talk at IWMS-16 in 2007 at Windsor, Ontario [see
the conference PowerPoint slides - see IWMS]. It became clear to us by 2008
that compounding 1EVs would produce larger 1EVs. Due to space restrictions
the compound work was not included in the LAA conference proceedings, with
the idea that it would be featured in a follow-up paper \cite{RCL} [RCL],
which now follows the present \ work.

So in the present context it is obvious that the compounded 1EV $n=4$\ Latin
and magic squares compound to $n=16$\ 1EV squares with $\mu =15$\ and rank $%
5 $, and thus worth raising to higher powers!

The idea is to place incremented versions of the basic square in positions
in a larger square with the same underlying pattern, which in this case is a
Latin square of 16 symbols which takes a lot of space which we abbreviate by
suggesting an order four filled with $\mathbf{A}=sud4a,\mathbf{B}=sud4a+4%
\mathbf{E}_{4},\mathbf{C}=sud4a+8\mathbf{E}_{4},\mathbf{D}=sud4a+12\mathbf{E}%
_{4}$:

\begin{equation}
\begin{bmatrix}
\mathbf{A} & \mathbf{B} & \mathbf{C} & \mathbf{D} \\ 
\mathbf{C} & \mathbf{D} & \mathbf{A} & \mathbf{B} \\ 
\mathbf{D} & \mathbf{C} & \mathbf{B} & \mathbf{A} \\ 
\mathbf{B} & \mathbf{A} & \mathbf{D} & \mathbf{C}%
\end{bmatrix}%
\end{equation}

or explicitly with blank rows and columns:

\begin{equation}
\begin{bmatrix}
1 & 2 & 3 & 4 &  & 5 & 6 & 7 & 8 &  & 9 & 10 & 11 & 12 &  & 13 & 14 & 15 & 16
\\ 
3 & 4 & 1 & 2 &  & 7 & 8 & 5 & 6 &  & 11 & 12 & 9 & 10 &  & 15 & 16 & 13 & 14
\\ 
4 & 3 & 2 & 1 &  & 8 & 7 & 6 & 5 &  & 12 & 11 & 10 & 9 &  & 16 & 15 & 14 & 13
\\ 
2 & 1 & 4 & 3 &  & 6 & 5 & 8 & 7 &  & 10 & 9 & 12 & 11 &  & 14 & 13 & 16 & 15
\\ 
&  &  &  &  &  &  &  &  &  &  &  &  &  &  &  &  &  &  \\ 
9 & 10 & 11 & 12 &  & 13 & 14 & 15 & 16 &  & 1 & 2 & 3 & 4 &  & 5 & 6 & 7 & 8
\\ 
11 & 12 & 9 & 10 &  & 15 & 16 & 13 & 14 &  & 3 & 4 & 1 & 2 &  & 7 & 8 & 5 & 6
\\ 
12 & 11 & 10 & 9 &  & 16 & 15 & 14 & 13 &  & 4 & 3 & 2 & 1 &  & 8 & 7 & 6 & 5
\\ 
10 & 9 & 12 & 11 &  & 14 & 13 & 16 & 15 &  & 2 & 1 & 4 & 3 &  & 6 & 5 & 8 & 7
\\ 
&  &  &  &  &  &  &  &  &  &  &  &  &  &  &  &  &  &  \\ 
13 & 14 & 15 & 16 &  & 9 & 10 & 11 & 12 &  & 5 & 6 & 7 & 8 &  & 1 & 2 & 3 & 4
\\ 
15 & 16 & 13 & 14 &  & 11 & 12 & 9 & 10 &  & 7 & 8 & 5 & 6 &  & 3 & 4 & 1 & 2
\\ 
16 & 15 & 14 & 13 &  & 12 & 11 & 10 & 9 &  & 8 & 7 & 6 & 5 &  & 4 & 3 & 2 & 1
\\ 
14 & 13 & 16 & 15 &  & 10 & 9 & 12 & 11 &  & 6 & 5 & 8 & 7 &  & 2 & 1 & 4 & 3
\\ 
&  &  &  &  &  &  &  &  &  &  &  &  &  &  &  &  &  &  \\ 
5 & 6 & 7 & 8 &  & 1 & 2 & 3 & 4 &  & 13 & 14 & 15 & 16 &  & 9 & 10 & 11 & 12
\\ 
7 & 8 & 5 & 6 &  & 3 & 4 & 1 & 2 &  & 15 & 16 & 13 & 14 &  & 11 & 12 & 9 & 10
\\ 
8 & 7 & 6 & 5 &  & 4 & 3 & 2 & 1 &  & 16 & 15 & 14 & 13 &  & 12 & 11 & 10 & 9
\\ 
6 & 5 & 8 & 7 &  & 2 & 1 & 4 & 3 &  & 14 & 13 & 16 & 15 &  & 10 & 9 & 12 & 11%
\end{bmatrix}%
\end{equation}

As expected from the compound spectral formulae of RCL this is diagonal
Latin [DL] with rank $5$. The results on powering-up are similar to $sud4a$
in Table 2, thus agreeing with our Proposition 2. but now with $C\%=55.8491$
and $Spread=1.76471$\ at $p=1$. Both $sud4a$ and its compound are DLs. Now
the linesums in each row are $(s_{16})^{p}=136^{p}$, i.e. for $p=1,2,3$
these are $136,18496,2515456.$ Their characteristic polynomials are
respectively: $x^{16}=136x^{15};x^{16}=18496x^{15};x^{16}=2,515,456x^{15}$,
etc. The first of these says that the eigenvalues are $136$ and fifteen
zeroes, i.e. another 1EV square, with multiplicity $\mu =15$, and so on.

Now we find an $n=16$ Jordan form, of which we show just an 8-by-8 portion
of the order 16 as the rest contains the only zero elements:

\begin{equation}
\mathbf{J}_{16}=%
\begin{bmatrix}
0 & 1 & 0 &  &  &  &  &  &  \\ 
0 & 0 & 1 &  &  &  &  &  &  \\ 
0 & 0 & 0 &  &  &  &  &  &  \\ 
&  &  & 0 & 1 & 0 &  &  &  \\ 
&  &  & 0 & 0 & 1 &  &  &  \\ 
&  &  & 0 & 0 & 0 &  &  &  \\ 
&  &  &  &  &  &  &  &  \\ 
&  &  &  &  &  &  &  &  \\ 
&  &  &  &  &  &  &  & 136%
\end{bmatrix}%
,
\end{equation}%
where two $3$-by-$3$ rank $2$\ Jordan blocks along the diagonal are followed
after nine rows and columns of zeros before the corner diagonal element of
the linesum, $136$. So we have a degree $2+2-1=3$ compound square, i.e.
raising it to the cube power will provide us with a constant matrix of the
previous form: $\mathbf{Z}^{k}=(\frac{l^{k}}{n})\mathbf{E}_{n}$, with the
observed rank decreasing by 2 each time because we have two Jordan blocks,
each contributing now a rank 1 decrease at each power. This Jordan structure
explains why constancy is attained at $p=3$, as for $sud4a$\ itself, and not
much higher as might have been expected by\textbf{\ }Cayley-Hamilton which
gives\textbf{\ }$\mathbf{A}^{16}=136\mathbf{A}^{15}$.

Continuing this compounding also gives larger Latin squares or orders $4\ast
16=64,4\ast 64=256,...$, all 1EV.

\section{n=16 compounding of associative 1EV $f360$}

Now we follow the previous compounding of $sud4a$\ with the analogous
approach for an $n=4$ 1EV magic square above. An example for a different
order four magic square is shown explicitly in CL. The Jordan form here
differs from the compound of $sud4a$ only by the magic linesum which here is 
$l_{16}=S_{16}=Sum(1..256)/16=2056.$

For the 1EV compounding of $f360$ the powers show constancy again at $p=3$
as for the previous $sud4a$ compounding, and the same reduction in rank as
for the compounding of $sud4a$, now with $C\%$ increasing from $64.3023$\
and Spread decreasing from $1.98444$\ at $p=1$ (compare with Table 6).$\ $

\subsection{Larger compound magic squares}

For 1EV compound magic squares combinations of the basic 1EV magic squares $%
n=4,5\;\&\;8$ featured in this work this means that compounding will produce
orders $4\ast 5=20$, $8\ast 5=40,16\ast 5=90$, etc., in addition to all
those that we could make just from powers of $n=4,5\;\&\;8$.

For $n_{1}=4$, $n_{2}=5$, the 1EV order $20$ compound magic squares have
rank $6$\ as expected [RCL], and power up to constancy at $p=4$, the same as
for the $n=5$\ 1EV component, as discussed earlier. This may be continued
for magic squares to multiples of $n=5$, and then products of powers of $4$
and $5$, \ i.e. $4\ast 5=20,16\ast 5=80,4\ast 25=100,..$.

\section{Products of different $n=4$\ 1EV magic squares}

Squaring a matrix is a pair product, \ and cubing a triple product. For the
associative Group III's we now try mixing different members of the same clan.

\subsection{Double and Triple products of different alpha squares}

[See Table 5 for details of the alpha and beta clans.]

Consider the following pair of alpha squares as a matrix product:

$pair1=f360.f790$; $\lambda _{i}=$ $\{1156,80,80,0\},d1=1316,d2=1156$, which
is no longer 1EV, but now 3EV and DA, which indicates that even with the
same clan the 1EV property is destroyed in this pair. We have tested powers
of $pair1$\ to 4th order and it is always 3EV, but slowly approaching
constancy.

However including either component $f360$ or $f790$\ as a third matrix
before or after $pair1$\ always renews the 1EV property:

$f360.pair1=f360.(f360.f790)$ 1EV, constant also at $p=2,3,...$

$pair1.f360=(f360.f790).f360$ 1EV, constant also at $p=3,4,...$

$pair1.f790=(f360.f790).f790$ 1EV, constant also at $p=2,3,...$

$f790.pair1=f790.(f360.f790)$ 1EV, constant also at $p=3,4,...$

These show that when the same square is adjacent the constancy appears in
the first (squared) power of the triple product, whereas when they are
separated by the other the constancy is delayed to the cube.

\subsection{Pair products of alpha and beta squares}

Multiplying a mixed pair of 1EV SVs, an $alpha$ square with a $beta$ square
also produces a 3EV $m=1$ square, e.g. try 1EV beta $f489$ with alpha $f790$%
: \ 

$pair2=f489.f790$; $\lambda _{i}=$ $\{1156,-32,16,0\},d1=1140,d2=1204$, so
no longer 1EV, but now 3EV DA, which again suggests that mixing clans
destroys the 1EV property.\bigskip\ 

Again we have tested powers of $pair2$\ to 4th order and it is always 3EV,
but slowly approaching constancy. However the following triples are now all
1EV, with the same pattern as for the previous cases which used a single
clan:

$f489.pair2=f489.(f489.f790)$ 1EV, constant also at $p=2,3,...$

$pair2.f489=(f489.f790).f489$ 1EV, constant also at $p=3,4,...$

$pair2.f790=(f489.f790).f790$ 1EV, constant also at $p=2,3,...$

$f790.pair2=f790.(f489.f790)$ 1EV, constant also at $p=3,4,...$

\subsection{Commutation}

Non-commuting (matrix) operators encountered in the matrix formulation of
quantum mechanics prompt us to consider this issue for the pair products
above. First with the same clan pair, $pair1=f360.f790$:

$comm1=f360.f790-f790.f360$: $\lambda _{i}=\{-80,80,0,0\},d1=0,d2=0$, so the
3EV $pair1$ has a 2EV commutator.

Second with the mixed clan pair, $pair2=bj489.f790$:

$comm2=bj489.f790-f790.bj489:\lambda _{i}=\{48,-32,-16,0\},$

$d1=0,d2=96$, and the\ 3EV $pair2$\ has a 3EV commutator.

With this insight we discontinue further exploration of these products.

\section{Fibonacci magic square ($\protect\mu =0$)}

We include this since CBH introduced Herta Freitag's order four magic square 
\cite{Freitag} even though it is not associative, nor are the elements in
arithmetic progression. Taking sixteen elements from the Fibonacci series,
and its $p=2$\ square:

\begin{equation}
\mathbf{A}=%
\begin{bmatrix}
13 & 89 & 97 & 34 \\ 
110 & 21 & 63 & 39 \\ 
68 & 94 & 55 & 16 \\ 
42 & 29 & 18 & 144%
\end{bmatrix}%
,\mathbf{A}^{2}=%
\begin{bmatrix}
17983 & 13130 & 12815 & 10361 \\ 
9662 & 17284 & 16160 & 11183 \\ 
15636 & 13660 & 15831 & 9162 \\ 
11008 & 10215 & 9483 & 23583%
\end{bmatrix}%
,...\text{.}
\end{equation}

These are non-singular we do not have a 1EV square:

$\lambda _{i}=$ $\{233,116.003,-68.1711,-47.8322\}$

The results for the first four powers are similar to Tables 3 and 7:

\begin{center}
$%
\begin{tabular}{|l|l|l|l|l|}
\hline
$p$ & $r$ & $C\%$ & spread & Type \\ \hline
$1$ & $4$ & $13.1301$ & $2.24893$ & magic \\ \hline
$2$ & $4$ & $38.8696$ & $1.06254$ & DA \\ \hline
$3$ & $4$ & $62.4492$ & $0.53231$ & DA \\ \hline
$4$ & $4$ & $78.655$ & $0.26527$ & DA \\ \hline
\end{tabular}%
\ \ \ \ \ $

Table 10 - Freitag Fibonacci square slowly converging towards constancy.
\end{center}

For the record we compounded this magic square to $n=16$, rank $7$, $%
C\%=52.2854$, $Spread=2.28586$.

\subsection{Prime 1EV Latin squares}

An example of a prime arithmetic progression is given on MathWorld \cite%
{prime}. From the ten terms given with a difference of $210$ we take the
first four for a prime Latin square modelled on $sud4a$ in (1) which is also
1EV,

\begin{equation}
\mathbf{A}=%
\begin{bmatrix}
199 & 409 & 619 & 829 \\ 
619 & 829 & 199 & 409 \\ 
829 & 619 & 409 & 199 \\ 
409 & 199 & 829 & 619%
\end{bmatrix}%
\text{.}
\end{equation}%
Here $\lambda _{i}=\{2056,0,0,0\};\sigma
_{i}=\{2056,840,420,0\},r=3,C\%=34.6522,Spread=1.22568$, and note the
closeness of these values to those for $sud4a$ which were $%
C\%=35.06,Spread=1.2$.

Clearly there are many more of these $n=4$ prime 1EV squares from the other
selections from the ten in MathWorld \cite{prime}, in addition to more from
compounding, and possibly for order $7,8,9,...$ Latin squares (we did not
find any for orders 5 and 6 in DMPS).

\section{Gerschgorin Disks for 1EV $n=4$ magic squares}

Following the discussion of the location of eigenvalues w.r.t to these disks
in Ortega \cite{Ortega}, we have not seen any previous use of these disks in
the literature of magic or Latin squares, and the opportunity to look at the
extreme case of just 1EV. \ It suffices to describe the image for $sud4a$
of\ four circles centered on the $x$-axis at $1,4,2,3$\ with radii $9,6,8,7$%
\ respectively (the element sums of the consecutive columns), eigenvalue $0$%
\ at the origin and eigenvalue $10$\ at the point where the four circles
touch.

\section{Conclusions\label{Conclusion}}

It was necessary to review previous "magic" terminology for these powers,
leading to our use of doubly-affine and diagonal doubly-affine to
characterize these magical squares.

Using doubly-affine in our \textbf{Proposition 1} to include both Latin and
magic squares we showed that certain associative $n=4,5$ 1EV squares power
up to constancy in a few steps. This was understood by using the
Cayley-Hamilton theorem for the corresponding characteristic equations,
leading to our \textbf{Proposition 2}.

CBH's Proposition 6.1\ had claimed that for \textbf{any} associative magic
squares their odd powers remained "magic", while their even powers were just
"semimagic", but we found exceptions. In addition we found evidence that
another type did satisfy their proposition. A different variation of CBH's
alternation between magic and semimagic found in some of our examples.

We found constancy for Latin $sud4a$ at $p=3$, as well as for its
compounding to $n=16$, \ and for some $n=4$ examples, while constancy
occured at $p=4$ for an $n=5$ associative 1EV magic square, as well as for
their compounding to multiplicative orders. It was necessary to study their
Jordan forms in order to understand constancy at $p=3$\ for their constancy
at $p=3$.

\section{Acknowledgements}

Miguel Amela drew our attention in September 2016, to the role of 1EV for
the even powers\ of our 8 order 4 associative magic squares by using the
transpose of $f360$. Walter Trump, our co-author in LAA for the 1eV squares.
Ron Nordgren \cite{Nordgren} for a copy of his recent work on Franklin
squares. Murray Alexander for suggesting that we look at the pair and triple
products.

[Author for correspondence: \textbf{email}: loly@umanitoba.ca]

\appendix

\end{document}